\documentclass[12pt,oneside,reqno,dvipsnames]{amsart}
\usepackage{graphicx}
\usepackage{caption,subcaption}
\usepackage{dcolumn}
\usepackage{amsmath,amsthm,amssymb}
\usepackage{booktabs} 
\renewcommand{\Re}{\operatorname{Re}}

\newcommand{\scaption}[1]{\caption{\small{#1}}}

\begin{document}

\title[Supercritical focusing wave equations outside a ball]{Dynamics
  at the threshold for blowup for supercritical wave equations outside
  a ball}

\author{Piotr Bizo\'n}
\address{Institute of Physics, Jagiellonian University, Krak\'ow,
  Poland}
\email{bizon@th.if.uj.edu.pl}

\author{Maciej Maliborski}
\address{Gravitational Physics, Faculty of Physics, University of
  Vienna, Boltzmanngasse 5, A-1090 Vienna, Austria}
\email{maciej.maliborski@univie.ac.at}

\thanks{This research was supported by the Polish National Science
  Centre grant no. DEC-2012/06/A/ST2/00397. Computations have been
  performed on Minerva cluster of the Max-Planck Institute for
  Gravitational Physics.}
\date{\today}
\begin{abstract}
  We consider spherically symmetric supercritical focusing wave
  equations outside a ball. Using mixed analytical and numerical
  methods, we show that the threshold for blowup is given by a
  codimension-one stable manifold of the unique static solution with
  exactly one unstable direction. We analyze in detail the convergence
  to this critical solution for initial data fine-tuned to the
  threshold.
\end{abstract}
\maketitle

\section{Introduction}
\label{sec:introduction}
This paper is concerned with the focusing semilinear wave equation for
a real scalar field $\phi(t,x)$,
\begin{equation}\label{eq-gen}
  \phi_{tt}=\Delta \phi + \phi^{2p+1}\,,
\end{equation}
outside a unit ball in $\mathbb{R}^d$ for odd $d\geq 3$. Here $p$ is a
positive integer greater than $\frac{2}{d-2}$ which corresponds to the
supercitical regime. We restrict ourselves to spherically symmetric
solutions $\phi(t,r)$, where $r=|x|$, satisfying the Dirichlet
boundary condition $\phi(t,1)=0$, hence we solve
\begin{equation}\label{eq}
\phi_{tt} = \phi_{rr}+\frac{d-1}{r} \phi_r + \phi^{2p+1} \quad \mbox{for}\,\, r\geq 1\,\,\mbox{with}\,\,\phi(t,1)=0,
\end{equation}
Initial data $(\phi(0,r), \phi_t(0,r))$ are assumed to be smooth and
compatible with the boundary condition.

Let us first briefly recall what is known about solutions of equation
\eqref{eq} in the whole space. For small initial data the solutions
are global in time and scatter to zero for $t\rightarrow \infty$
\cite{st}. The behavior of large solutions is only partially
understood. In the case $d=3$, the numerical studies reported in
\cite{bct} show that for generic large initial data the solutions blow
up as $u \sim (T-t)^{-1/p}$ for $t\nearrow T<\infty$. The nonlinear
stability of this ODE blowup was proved by Donninger \cite{d}. In
addition, there exists a countable family of unstable self-similar
solutions which correspond to non-generic finite time blowups
\cite{bmw} and the unique self-similar solution with exactly one
unstable direction was shown numerically to be critical in the sense
that its codimension-one stable manifold separates dispersive and
singular solutions \cite{bct}. The codimension-one nonlinear stability
of this critical solution was proved by Donninger and Sch\"orkhuber
\cite{ds} (see also \cite{gs} for an analogous result in higher
dimensions).

The presence of the obstacle does not affect the qualitative behavior
of generic solutions, that is small solutions scatter to zero, while
large solutions exhibit the ODE blowup. However, the obstacle breaks
the scaling symmetry thereby excluding self-similar solutions and at
the same time allowing for static solutions. These static solutions
are known from studies of elliptic equations \cite{c} but, as far as
we know, their role in dynamics has not been studied\footnote{Note
  added: while completing this paper, we were informed by Thomas
  Duyckaerts about his work with J. Yang in which they proved that any
  global-in-time solution of equation \eqref{eq} either scatters to
  zero or converges (up to a dispersive term) to one of the static
  solutions. The proof is based on the concentration-compactness
  technique which gives no information about the rate of
  convergence.}. For completeness, in the next section we give an
elementary proof of existence of a countable family of static
solutions with increasing number of nodes. We also prove that the
nodal index of these solutions counts the number of their unstable
modes. The main goal of this paper is to show, using mixed numerical
and analytical methods, that the static solution with one unstable
mode plays the role a critical solution whose codimension-one stable
manifold separates dispersive and singular solutions.

\section{Static solutions and their stability}
\label{sec:stat-solut-their}
\begin{figure}[h!]
  \centering
  \includegraphics[width=0.7\textwidth]{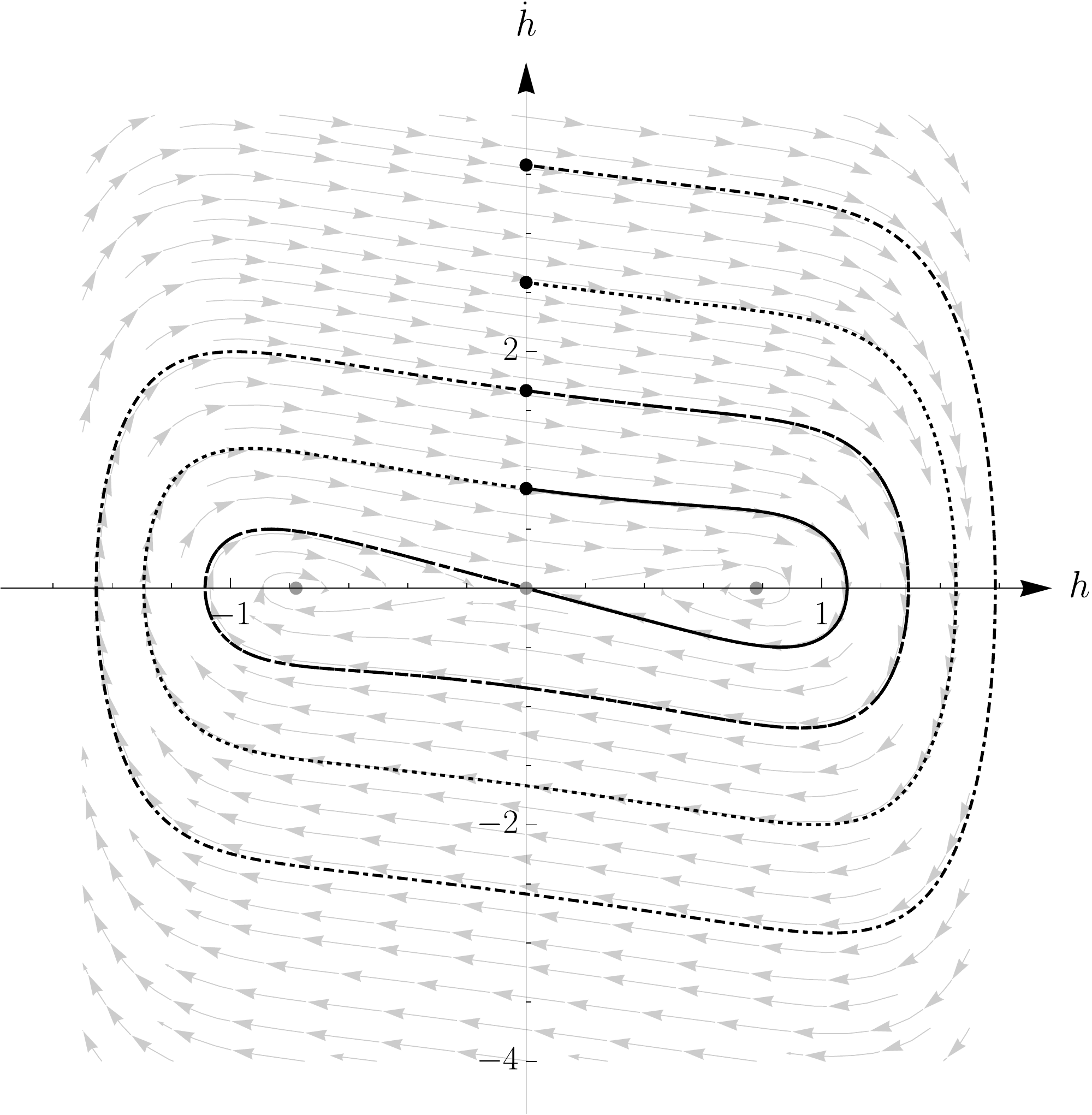}
  \scaption{Phase portrait $(h,\dot{h})$ for a sample pair
    $(d,p)=(3,3)$. The trajectories of the first four static solutions
    $h_{n}$ ($n=0,1,2,3$) are plotted with distinct line styles and
    their starting points $(h_{n}(0),\dot{h}_{n}(0))=(0,b_{n})$ are
    marked with black circles.}
  \label{fig:1}
\end{figure}

For time-independent solutions equation \eqref{eq} reduces to the
radial Lane-Emden equation
\begin{equation}\label{eqs}
  \phi_{tt}=\phi_{rr}+\frac{d-1}{r}\, \phi_r + \phi^{2p+1},
\end{equation}
which after the change of variables (introduced by Fowler in \cite{F})
\begin{equation}\label{sh}
  s=\ln{r},\quad h(s)=r^{1/p} \phi(r)
\end{equation}
transforms into the autonomous ordinary differential equation
($\cdot=d/ds$)
\begin{equation}\label{eqh}
 \ddot h+\left(d-2-\frac{2}{p}\right)  \dot h - \frac{1}{p} \left(d-2-\frac{1}{p}\right) h + h^{2p+1}=0\,.
\end{equation}
For $p>2/(d-2)$ the `friction' coefficient in \eqref{eqh} is positive
and from an elementary phase-plane analysis (see the phase portrait in
Fig.~\ref{fig:1}) it follows that there exist infinitely many initial
values $(h_n(0),\dot h_n(0))=(0,b_n)$, where $n$ is a nonnegative
integer, for which the phase trajectory makes $(n+1)$ half rotations
around the origin and then tends to the saddle point at the origin
along the stable direction
\begin{equation}\label{finf}
  h_n(s) \sim c_n e^{{-(d-2-1/p)s}}\quad \mbox{for}\,\, s\rightarrow \infty.
\end{equation}
In terms of the original variables these trajectories correspond to
finite energy static solutions $\phi_n(r)$ which vanish at $r=1$ and
decay as $c_n/r^{d-2}$ for $r\rightarrow \infty$. The first few values
of parameters $b_n$ and $c_n$ determined numerically for several pairs
$(d,p)$ are given in Tab.~\ref{tab:1}.

\vskip 0.2cm
\begin{table}
  \begin{tabular}{|c|ccc|}
    \toprule
    $(d,p)$ & $(b_{0},c_{0})$ & $(b_{1},c_{1})$ & $(b_{2},c_{2})$ \\
    \midrule
    (3,3) & (0.84261,\,-4.46847) & (1.67035,\,\phantom{-}21.7658) & (2.58523,\,-62.5081) \\
    (3,4) & (1.20653,\,-3.71646) & (2.48958,\,\phantom{-}13.0365) & (3.90145,\,-28.9009) \\
    (3,5) & (1.41849,\,-3.35818) & (2.95061,\,\phantom{-}10.1979) & (4.61581,\,-20.3151) \\
    (5,1) & (5.51059,\,-22.5426) & (12.4733,\,\phantom{-}209.872) & (21.5494,\,-1005.52) \\
    (5,2) & (7.70805,\,-8.22701) & (18.1434,\,\phantom{-}32.8788) & (30.9438,\,-79.2027) \\
    (5,3) & (7.69629,\,-5.64440) & (17.4958,\,\phantom{-}17.8598) & (28.8616,\,-36.3276) \\
    \bottomrule
  \end{tabular}
  \vskip 1ex
  \scaption{Parameters $b_{n}$ and $c_n$ of the first three static
    solutions for a few pairs $(d,p)$.}
  \label{tab:1}
\end{table}

We remark that no static solutions exist in the critical and
subcritical cases $p\leq 2/(d-2)$, as follows, for instance, from the
identity
\begin{equation}\label{poh}
  \frac{1}{2} \dot h^2(0)= \left(d-2-\frac{2}{p}\right) \int_0^{\infty} \dot h^2(s) ds\,,
\end{equation}
which arises from multiplying equation \eqref{eqh} by $\dot h$ and
integrating by parts.

\vskip 0.2cm

The role of static solutions in dynamics depends on their stability
properties. To determine the linear stability of the solution
$\phi_n(r)$, we substitute $\phi(t,r)=\phi_n(r)+w(t,r)$ into
\eqref{eq}. Dropping nonlinear terms in $w$, we get the linearized
equation
\begin{equation}\label{eq-w}
  w_{tt} = w_{rr}+\frac{d-1}{r} w_r + (2p+1) \phi_n^{2p} w\,.
\end{equation}
Substituting $w(t,r)=e^{\lambda t} v(r)$ into \eqref{eq-w}, we obtain
the eigenvalue problem
\begin{equation}\label{eigen-v}
  L_n v := \left(-\frac{d^2}{dr^2}-\frac{d-1}{r} \frac{d}{dr} -(2p+1) \phi^{2p}_n(r)\right) v = -\lambda^2 v.
\end{equation}
For each $n$ the operator $L_n$ is essentially self-adjoint in the
Hilbert space
$X=\{v: \int_1^{\infty} v^2(r) r^{d-1} dr<\infty, \, v(1)=0\}$. Since
$\phi_n(r)$ is bounded and decays to zero at infinity, $L_{n}$ has a
continuous spectrum $[0,\infty)$. Note that the function generated by
scaling
\begin{equation} \label{zeromode}
  v_0^{(n)}(r)=\Big.\frac{d}{d\alpha} \alpha^{1/p} \phi_n(\alpha r)\Big|_{\alpha=1}=r\phi_n'(r) + \frac{1}{p} \phi_n(r)
\end{equation}
solves equation \eqref{eigen-v} for $\lambda=0$ (but it is not an
eigenfunction because it does not belong to $X$).  From the
phase-plane analysis above it follows that $v_0^{(n)}(r)$ has exactly
$n$ zeros which implies by the Sturm oscillation argument that the
operator $L_n$ has exactly $n$ negative eigenvalues, hereafter denoted
by $-(\lambda_k^{(n)})^2$ ($k=1,..,n$). Consequently, the static
solution $\phi_n$ has exactly $n$ unstable modes
$v_k^{(n)}(r) e^{\lambda_k^{(n)} t}$.
In what follows we focus on dynamics near the ground state solution
$\phi_0(r)$ which has exactly one unstable mode
$v_1(r) e^{\lambda_1 t}$ (henceforth we drop the superscript $(n)$ on
the eigenvalues and eigenfunctions).
Due to the presence of the unstable mode, generic solutions of the
linearized equation \eqref{eq-w} grow exponentially. This instability
can be eliminated by preparing initial data that are orthogonal to the
unstable mode. The solutions starting from such special initial data
decay in time due to a combination of two dispersive effects: the
quasinormal ringdown and the polynomial tail. The rate of decay of the
tail is determined by the fall-off of the potential term in
\eqref{eq-w}: since $\phi_0^{2p}(r)\sim r^{-2p(d-2)}$ for
$r\rightarrow \infty$, it follows that $\phi(t,r) \sim t^{-\beta}$,
where $\beta=d-4+2p(d-2)$, for any fixed $r>1$ and
$t\rightarrow \infty$ \cite{st,bcr}. The ringdown is determined by the
quasinormal modes which are solutions of the eigenvalue equation
\eqref{eigen-v} with $\Re(\lambda) < 0$ satisfying the outgoing wave
condition $v(r) \sim e^{-\lambda r}$ for $r\rightarrow \infty$.  As
the concept of quasinormal modes is inherently related to the loss of
energy by radiation, the unitary evolution \eqref{eq-w} and the
associated self-adjoint eigenvalue problem \eqref{eigen-v} do not
provide a natural setting for analysing quasinormal modes, both from
the conceptual and computational viewpoints. For this reason we
postpone the discussion of quasinormal modes until the next section
where a new nonunitary formulation will be introduced.

\section{Characteristic initial-boundary value formulation}
\label{sec:char-init-bound}
The rest of the paper is devoted to dynamics of convergence to
$\phi_0$ for initial data fine-tuned to the threshold. To this order
we introduce the null coordinate $u=t-r$ and the inverse radial
coordinate $x=1/r$ which compactifies the spatial domain to the
interval $0\leq x \leq 1$. Then $f(u,x)=r^{\frac{d-1}{2}} \phi(t,r)$
satisfies the equation
\begin{equation}\label{eqf}
  2 f_{ux} +x^2 f_{xx}+2x f_x -\frac{1}{4}(d-3)(d-1) f + x^{\alpha} f^{2p+1}=0,\quad f(u,1)=0,
\end{equation}
where $\alpha=\frac{1}{2}(p(d-1)-2)$. We note in passing that
equation \eqref{eqf} can be written as the conservation law
\begin{equation}\label{claw}
  \partial_x\left(f_u^2+x^2 f_u f_x\right)=\partial_u\left(\frac{1}{2} x^2 f_x^2+\frac{1}{8}(d-3)(d-1) f^2-\frac{1}{2p+2} x^{\alpha} f^{2p+2}\right)\,,
\end{equation}
which upon integration gives the energy loss formula
\begin{equation}\label{bondi}
  \frac{dE}{du}=-f_u(u,0)^2,
   \end{equation}
   where
   \begin{equation}\label{energy}
   E[f]=\int_0^1 \left(\frac{1}{2} x^2 f_x^2 +\frac{1}{8} (d-3)(d-1) f^2-\frac{1}{2p+2} x^{\alpha} f^{2p+2} \right) dx\,.
\end{equation}
As follows from section~\ref{sec:stat-solut-their}, equation
\eqref{eqf} has infinitely many static solutions $f_n(x)$ which behave
as $f_n(x)\sim c_n x^{\frac{d-3}{2}}$ near $x=0$ and vanish at $x=1$.
These static solutions are critical points of the energy functional
$E[f]$; in particular, the solution $f_0$ is the ground state.

We now repeat the linear stability from the previous section by
substituting $f(u,x)=f_0(x)+e^{\lambda u} v(x)$ into \eqref{eqf} and
linearizing. This yields the eigenvalue problem
\begin{equation}\label{eqv}
  x^2 v''+2x v'+2\lambda v' -\frac{1}{4} (d-3)(d-1) v + (2p+1) x^{\alpha} f_{0}^{2p}(x) v=0, \quad \, v(1)=0.
\end{equation}
An advantage of this formulation is that it allows us to treat
quasinormal modes as genuine eigenfunctions. To do so we must specify
the desired behavior of eigenfunctions at $x=0$ which is rather subtle
because this endpoint is an essential singularity. The two linearly
independent solutions of equation \eqref{eqv} near $x=0$ have the
following leading behaviors
\begin{equation}\label{eqv-asympt}
  v_{g}(x) \sim 1, \qquad v_{b}(x) \sim e^{2\lambda/x}\,,
\end{equation}
where the subscripts $g$ and $b$ stand for `good' and `bad' solutions,
respectively. At $x=0$ the solution $v_{g}(x)$ admits a formal Taylor
series, while the solution $v_{b}(x)$ has an essential singularity. In
terms of the original variables, these two solutions correspond to the
outgoing and ingoing waves, respectively, thus we demand that the
eigenfunctions have no admixture of $v_{b}$. Having a good solution
near $x=0$, one can shoot it towards $x=1$ and determine the
eigenvalues from the boundary condition $v(1)=0$.  Since the formal
Taylor series of $v_{g}$ is in general divergent, in practice we take
the asymptotic expansion of $v_g(x)$ at some small $x_0$ and truncate
it at the least term. While this optimal truncation approach works
very well in the case of positive (unstable) eigenvalues, it is not
precise enough for the eigenvalues with $\Re(\lambda)<0$ because in
this case the bad solution $v_{b}(x)$ is smaller than any power of $x$
for $x\rightarrow 0^+$. To capture such a small term we use the Borel
summation method which goes as follows.
Given a formal power series $v_{g}(x)= \sum_{k=0}^{\infty}a_{k}x^{k}$,
we Borel transform it
\begin{equation}
  \label{eq:8}
  \mathcal{B}(x)=\sum_{k=0}^{\infty}\frac{a_{k}}{k!}x^{k}\,.
\end{equation}
and then take the Laplace transform to get the Borel sum
\begin{equation}
  \label{eq:9}
  \mathcal{B}_{S}(x) = \int_{0}^{+\infty}e^{-t}\mathcal{B}(tx) dt\,.
\end{equation}
In practice, we truncate the series in (\ref{eq:8}) at some high order
$2K$ and accelerate the convergence by using the (diagonal) Pad\'e
approximation $P_{K}\mathcal{B}$. Because of possible poles of the
Pad\'e approximation on the real axis, we deform the integration
contour in \eqref{eq:9} by introducing the following path on the
complex plane $\gamma(\varepsilon)=\gamma_{1}\cup\gamma_{2}$ with
\begin{equation}
  \label{eq:10}
  \gamma_{1}\,:\ [0,\varepsilon]\ni s\rightarrow is\in\mathbb{C},
  \quad
  \gamma_{2}\,:\ [0,\infty)\ni s\rightarrow i\varepsilon+s\in\mathbb{C},
\end{equation}
where $\varepsilon$ is a free real parameter\footnote{In practice,
  having an initial guess for the eigenvalue $\lambda$ (based on the
  optimal truncation method), we looked at the distribution of poles
  of $P_{K}\mathcal{B}(z)$ on the complex $z$-plane to estimate the
  value of the parameter $\varepsilon$. In most cases $\varepsilon=1$
  worked reasonably well.} (since the integrand decays sufficiently
fast we need not to close the contour `at infinity'). In our
calculations we took the mean of integrals along the contours
$\gamma(\varepsilon)$ and $\gamma(-\varepsilon)$, thus we approximated
(\ref{eq:9}) by
\begin{equation}
  \label{eq:11}
  \mathcal{B}_{S}(x) \approx \frac{1}{2}\int_{\gamma(\varepsilon)}
  e^{-t}P_{K}\mathcal{B}(tx) dt
  + \frac{1}{2}\int_{\gamma(-\varepsilon)}
  e^{-t}P_{K}\mathcal{B}(tx) dr\,.
\end{equation}
and then computed these integrals numerically.

Having set up initial conditions at $x_0$ (either by the optimal
truncation or Borel summation), we integrated Eq.~(\ref{eqv}) using an
adaptive Runge-Kutta method of 8th order and then determined the
eigenvalues by solving the boundary condition $v(1)=0$ with Newton's
method (see Tab.~\ref{tab:2}). To suppress round-off errors we used an
extended precision arithmetics, typically with more than $20$
digits. Particularly demanding was the computation of the first stable
eigenvalue $\lambda_2$ for $d=3$. For example, in the case
$(d,p)=(3,3)$ we used $x_{0}\approx 0.00825$, $K=128$,
$\varepsilon = 10$ and Gauss quadratures with $38$ and $128$ nodes for
Gauss-Legendre and Gauss-Laguerre rules to compute the integrals
(\ref{eq:11}) along $\gamma_{1}$ and $\gamma_{2}$ respectively. This
scheme provided an accurate enough initial conditions for the shooting
algorithm to produce $\lambda_2\approx -0.04328358...$ whose first 15
digits did not depend on the choice of the starting point $x_{0}$
which made us feel confident that the result is correct.

\begin{table}[t!]
  \centering
  \begin{tabular}{|c|ccc|}
    \toprule
    $(d,p)$ & $\lambda_{1}$ & $\lambda_{2}$ & $\lambda_{3}$ \\
    \midrule
    (3,3) & 0.4376132 & -0.04328358 & -0.7359469 $\pm$ 0.6611351$\,i$ \\
    (3,4) & 0.9119156 & -0.12566311 & -0.9112554 $\pm$ 1.228442$\,i$\phantom{0} \\
    (3,5) & 1.393964\phantom{0}  & -0.21578421 & -0.9589717 $\pm$ 1.608909$\,i$\phantom{0} \\
    (5,1) & 1.412962\phantom{0} & -0.1580264 $\pm$ 0.2094073$\,i$ & -3.663357\phantom{0} $\pm$ 1.863078$\,i$\phantom{0} \\
    (5,2) & 4.006646\phantom{0} & -0.5943277 $\pm$ 0.4789266$\,i$ & -5.062170\phantom{0} $\pm$ 5.850155$\,i$\phantom{0} \\
    (5,3) & 6.472988\phantom{0} & -0.9450331 $\pm$ 0.5032462$\,i$ & -5.050332\phantom{0} $\pm$ 8.049461$\,i$\phantom{0} \\
    \bottomrule
  \end{tabular}
  \vskip 1ex
  \scaption{The unstable eigenvalue $\lambda_1$ and two least damped
    stable eigenvalues $\lambda_2, \lambda_3$ of the linearized
    operator around the ground state solution $f_{0 }$ for several
    pairs $(d,p)$. All given digits are significant.}
  \label{tab:2}
  \vspace{-0.5cm}
\end{table}

In the appendix we describe a different method of finding the spectrum of the linearized problem which reproduces all the above eigenvalues except for those that lie on the negative real axis.
\section{Critical evolution}
\label{sec:critical-evolution}
In this section we give numerical evidence supporting our conjecture that the ground state solution $f_0$ sits at the threshold for generic blowup.

Before presenting results we briefly describe our method of solving
numerically the initial-boundary value problem \eqref{eqf}. We use the
method of lines with a spectral element method for space
discretization.  The starting point of this approach is a weak
formulation of equation (\ref{eqf}). The spatial domain is divided
into non-overlapping intervals and on each interval the integral is
approximated using the Gauss-Legendre quadrature formula.  We
typically use 16 grid points in each of 9 equal size intervals of the
spatial domain. The coupling between the intervals is enforced by the
requirement of smoothness. At the sphere $x=1$ we impose the Dirichlet
condition, while at null infinity $x=0$ no condition is imposed. The
resulting equations are integrated in time using the 6th order
Runge-Kutta scheme with a fixed time step.
The presence of the mixed derivative in equation (\ref{eqf}) required
a solution of the algebraic system at the internal steps of the
Runge-Kutta scheme.

To get a clear picture of near critical evolution it was instrumental
to use high precision arithmetics which is computationally
expensive. The numerical algorithm described above gave satisfactory
results at an acceptable cost. The efficiency of the spectral element
method is due to its fast convergence and the sparse (block diagonal)
structure of matrices.  To further speed up calculations we use a
parallel version of the bisection search to fine tune the initial
data. The code was written in Mathematica.

\vskip 0.2cm

\begin{figure}[h!]
  \centering
  \includegraphics[width=0.95\textwidth]{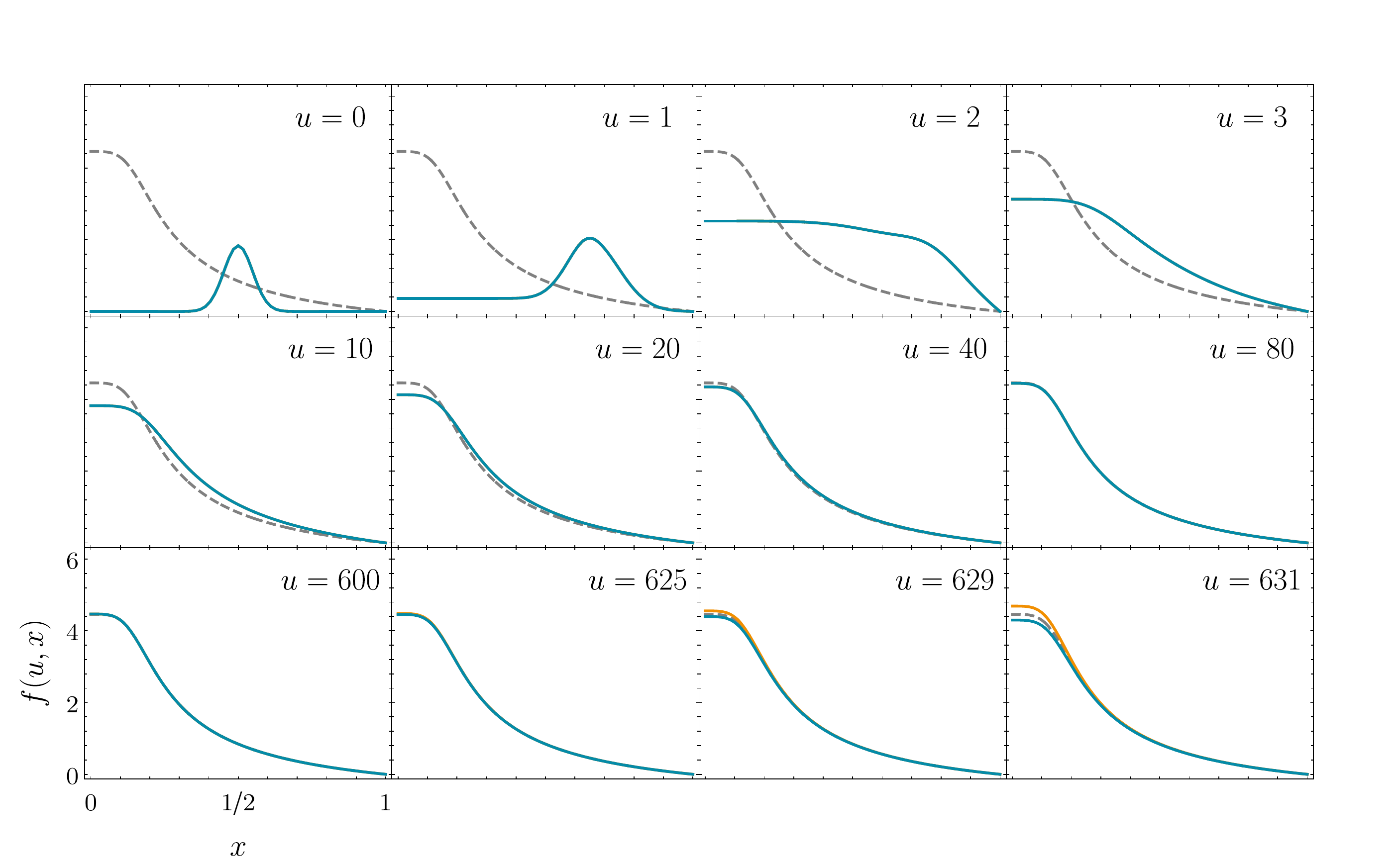}
  \vspace{2ex}
  \includegraphics[width=0.95\textwidth]{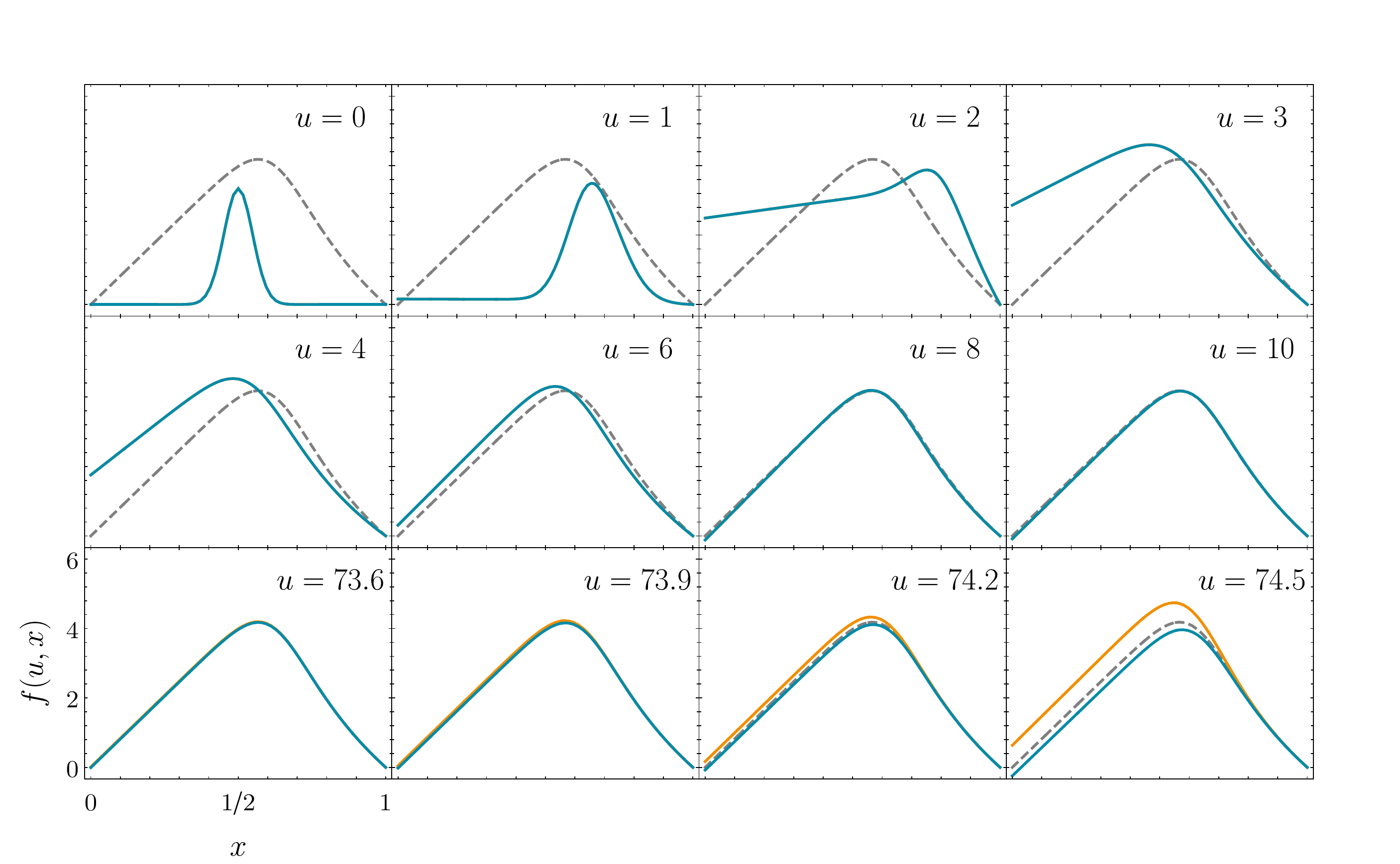}
  \scaption{Snapshots of near critical evolution of initial data
    \eqref{eq:12} for $(d,p)=(3,3)$ (upper plot) and $(d,p)=(5,2)$
    (lower plot). The amplitudes of marginally subcritical (blue
    lines) and supercritical (orange lines) data differ by
    $10^{-128}$. Initially the solutions evolve together, approach
    $f_0$ (dashed lines) for intermediate times, and eventually depart
    in opposite directions.}
  \label{fig:10}
\end{figure}

We illustrate our numerical results for a one-parameter family of
initial data
\begin{equation}
  \label{eq:12}
  f(0,x) = a\sin^{2}\left(\pi x\right)e^{-200(x-1/2)^{2}}\,.
\end{equation}
which interpolates between dispersion to zero for small amplitudes $a$
and the ODE blowup for large $a$. Using bisection we fine tune the
amplitude to the critical value $a_*$ separating these two generic
behaviors. For such fine-tuned initial data we observe for
intermediate times the convergence to the ground state $f_0$. This is
shown in Fig.~\ref{fig:10} for two pairs $(d,p)=(3,3)$ and $(5,2)$.

\begin{figure}[h]
  \centering
  \includegraphics[width=0.48\textwidth]{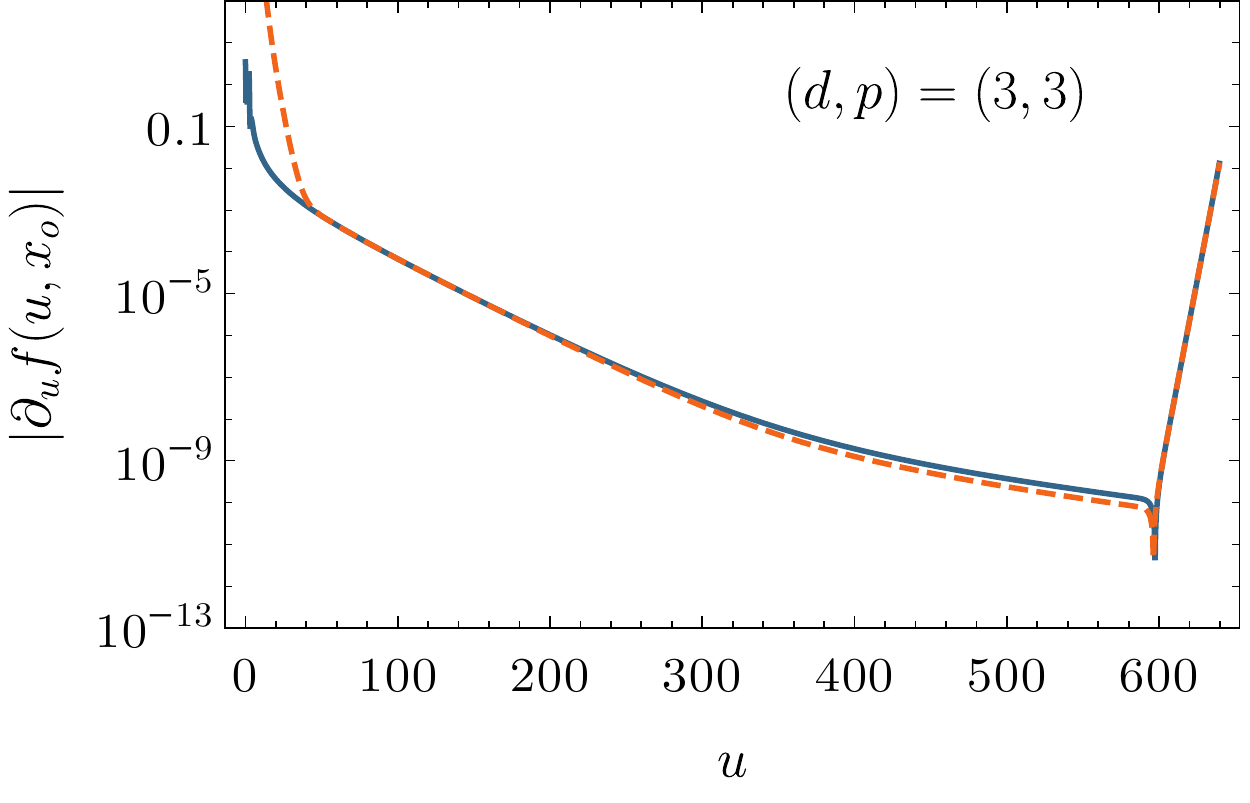}
  \hspace{1.5ex}
  \includegraphics[width=0.48\textwidth]{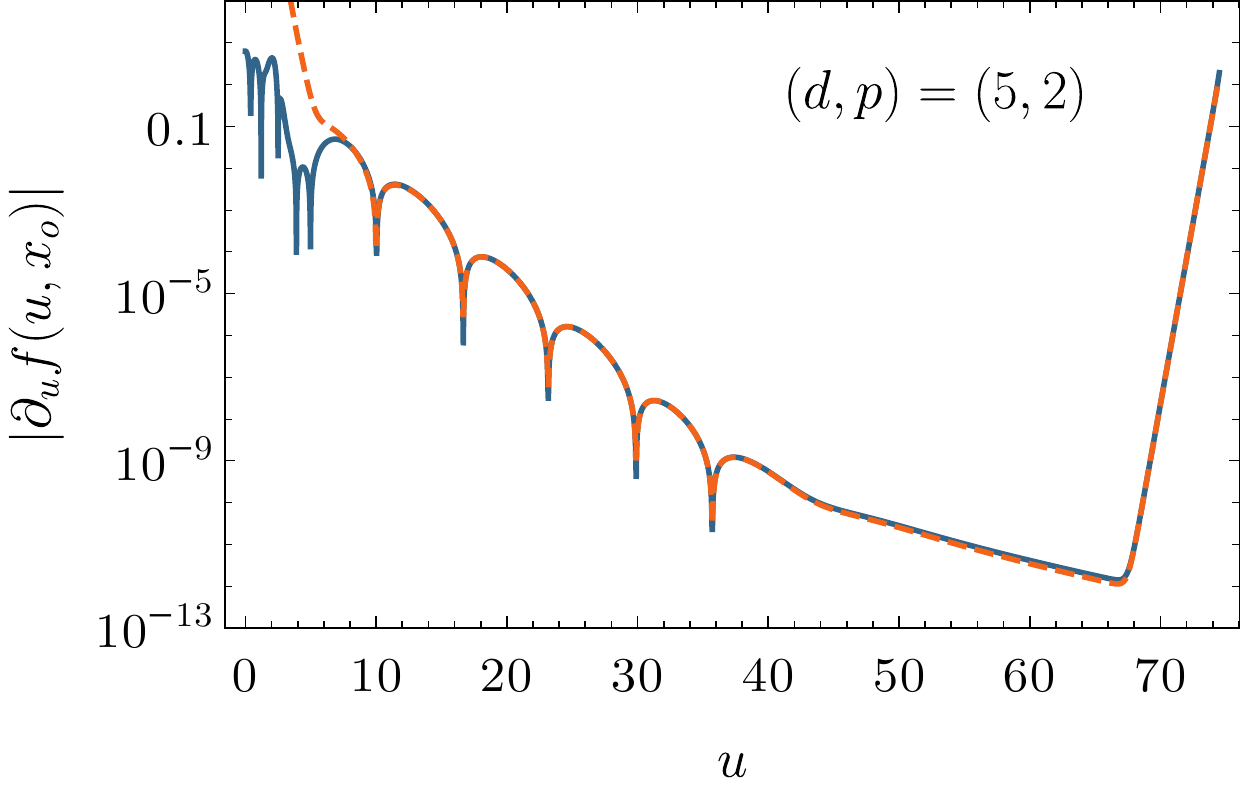}
  \scaption{Pointwise convergence to $f_0$ at a sample interior point
    $x_{o}=5/9$ for the same marginally subcritical data as in
    Fig.~\ref{fig:10}. The dashed lines depict the fits based on
    formula (\ref{eq:13}). For $(d,p)=(5,2)$ the evolution has two
    well-separated phases: the quasinormal ringdown followed by the
    polynomial tail. For $(d,p)=(3,3)$ the least damped quasinormal
    mode is non-oscillatory and very slowly decaying which makes it
    harder to separate it from the tail.}
  \label{fig:2}
\end{figure}

For intermediate times, when the nearly critical solution is close to
the ground state $f_{0}$, the dynamics is well approximated (for any
fixed $x$) by the linearized formula
\begin{equation}
  \label{eq:13}
  f(u,x) = f_{0}(x) + c_{1}e^{\lambda_{1}u} +
  \underbrace{\Re\left(c_{2}e^{\lambda_{2}u}\right) + \ldots}_\text{ringdown}
  + \underbrace{c_{3}u^{-\beta}\left(1 + \frac{c_{4}}{u} +
      \frac{c_{5}}{u^2}+\ldots\right)}_\text{tail}
  \,,
\end{equation}
where $c_{i}$ are ($x$ dependent) parameters and dots denote
subleading terms. For exactly critical data the coefficient $c_{1}$
vanishes. Our bisection procedure ensures that $c_{1}\sim a-a_{*}$ is
very small, typically of order $ 10^{-128}$, which gives a reasonably
long span of time $\propto -\frac{1}{\lambda_1} \log|a-a_*|$ over
which the linearized approximation \eqref{eq:13} is expected to hold
and can be fitted to the nearly critical solution shown in
Fig.~\ref{fig:10}. Performing this fit (keeping the exponent of the
tail $\beta=d-4+2p(d-2)$ fixed) we reproduce the eigenvalues
$\lambda_1$ and $\lambda_2$ with precision of $0.01\%$ which is very
reassuring; see Fig.~\ref{fig:2}.

\appendix

\section{\small{Pseudospectral solution of the linear problem}}
\label{sec:pseud-solut-line}
Here we present a simple algebraic method of solving the linearized
characteristic initial-boundary value problem which reproduces most
(but not all) results from section~\ref{sec:char-init-bound}.

\begin{figure}[h]
  \centering
  \includegraphics[width=0.55\textwidth]{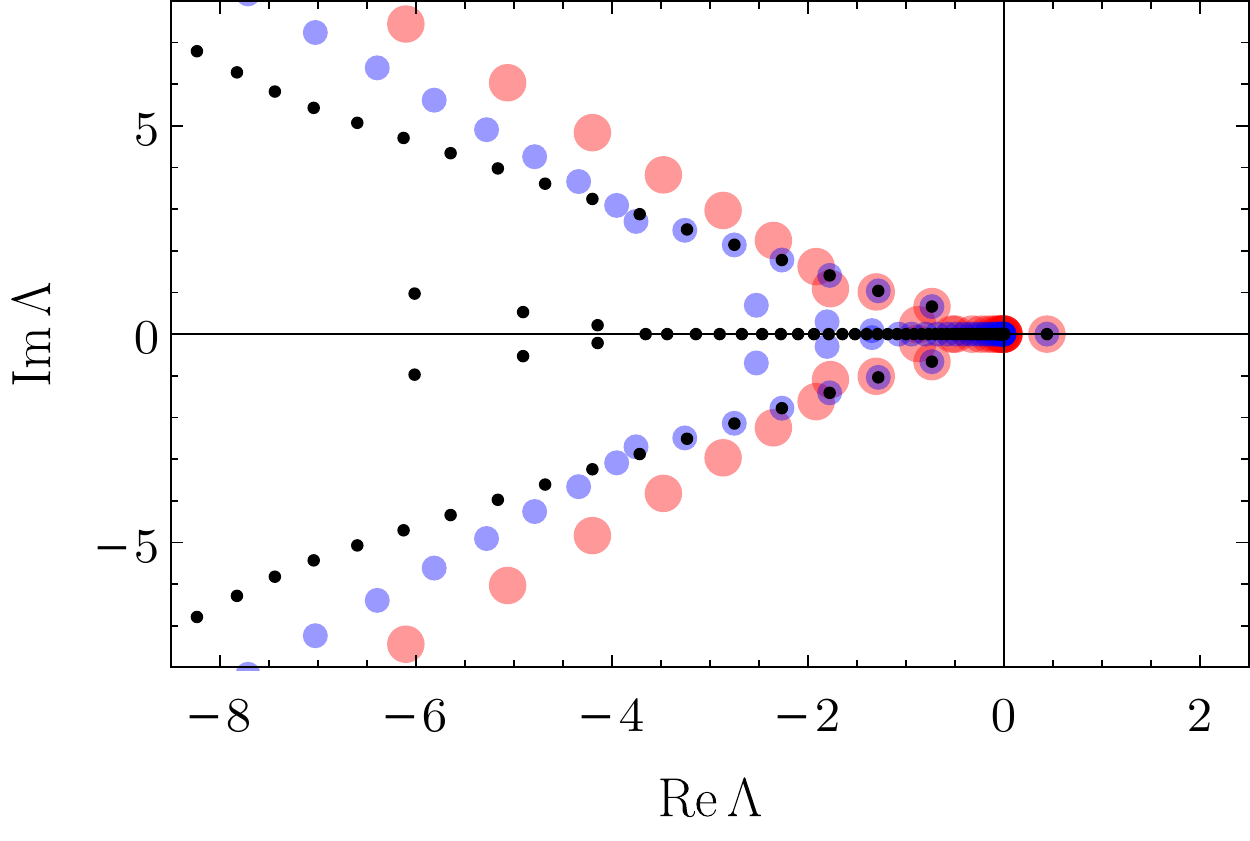}
  \scaption{Eigenvalues $\Lambda_i$ of the discrete operator
    $\mathbf{D}^{-1}\mathbf{L}_{0}$ in the $(d,p)=(3,3)$ case for
    different numbers of Chebyshev polynomials used in the
    approximation: $N=64$ (red), $N=128$ (blue), $N=256$ (black). As
    $N$ increases, we observe accumulation of spurious eigenvalues on
    the negative real axis and convergence to the genuine eigenvalues
    elsewhere.}
  \label{fig:a1}
\end{figure}

Linearization of equation \eqref{eqf} around a static solution $f_{n}$
yields
\begin{equation}
  \label{eq:1}
  \partial^{2}_{ux}v = L_{n}v\,,
\end{equation}
with
\begin{equation}
  \label{eq:2}
  L_{n}:=\frac{1}{2}\left(-x^{2}\partial_{x}^{2}-2x\partial_{x} +
    \frac{1}{4}(d-3)(d-1) - (2p+1) x^{\alpha}f_{n}^{2p}(x) \right)\,.
\end{equation}
Discretization in space transforms equation (\ref{eq:1}) into a system
of $N$ coupled constant-coefficient ODEs, where $N$ is the number of
degrees of freedom introduced by discretization. In the case at hand,
$N$ is the number of Chebyshev polynomials used in the spatial
approximation of $v(u,x)$. This semi-discrete problem has the form
\begin{equation}
  \label{eq:3}
  \mathbf{D}\frac{d}{du}\mathbf{v} =
  \mathbf{L}_{n}\mathbf{v}\,,
\end{equation}
where $\mathbf{v}$ is a vector of unknowns, $\mathbf{D}$ is an
invertible discrete version of $\partial_{x}$ which incorporates the
boundary condition $v(u,1)=0$, and $\mathbf{L}_{n}$ is a
discretization of $L_n$. We rewrite (\ref{eq:3}) as
\begin{equation}
  \label{eq:4}
  \frac{d}{du}\mathbf{v} = \mathbf{D}^{-1}\mathbf{L}_{n}\mathbf{v}\,,
\end{equation}
and performing diagonalization
\begin{equation}
  \label{eq:6}
  \mathbf{D}^{-1}\mathbf{L}_{n} = \mathbf{P}\Lambda\mathbf{P}^{-1}, \quad
  \Lambda = \mathrm{diag}\left(\Lambda_{1},\ldots,\Lambda_{N}\right)\,,
\end{equation}
we solve the system (\ref{eq:3}) by exponentiation
\begin{equation}
  \label{eq:5}
  \mathbf{v}(u) = \mathbf{P}\exp\left(
    \Lambda u\right)\mathbf{P}^{-1}\mathbf{v}(0),
\end{equation}
where $\mathbf{v}(0)$ is a vector of initial data. As $N$ grows, the
eigenvalues $\Lambda_{i}$, $i=1,\ldots,N$ tend to the eigenvalues
$\lambda$ of (\ref{eqv}), hence by increasing $N$ we uncover more and
more eigenvalues found in section~\ref{sec:char-init-bound} with the
shooting method. This is illustrated in Fig.~\ref{fig:a1}. The
drawback of this method is the accumulation of spurious eigenvalues on
the negative real axis which makes it hardly possible to extract the
genuine eigenvalues lying on that axis (such as $\lambda_2$ for
$d=3$).

\begin{figure}[h]
  \centering
  \includegraphics[width=0.48\textwidth]{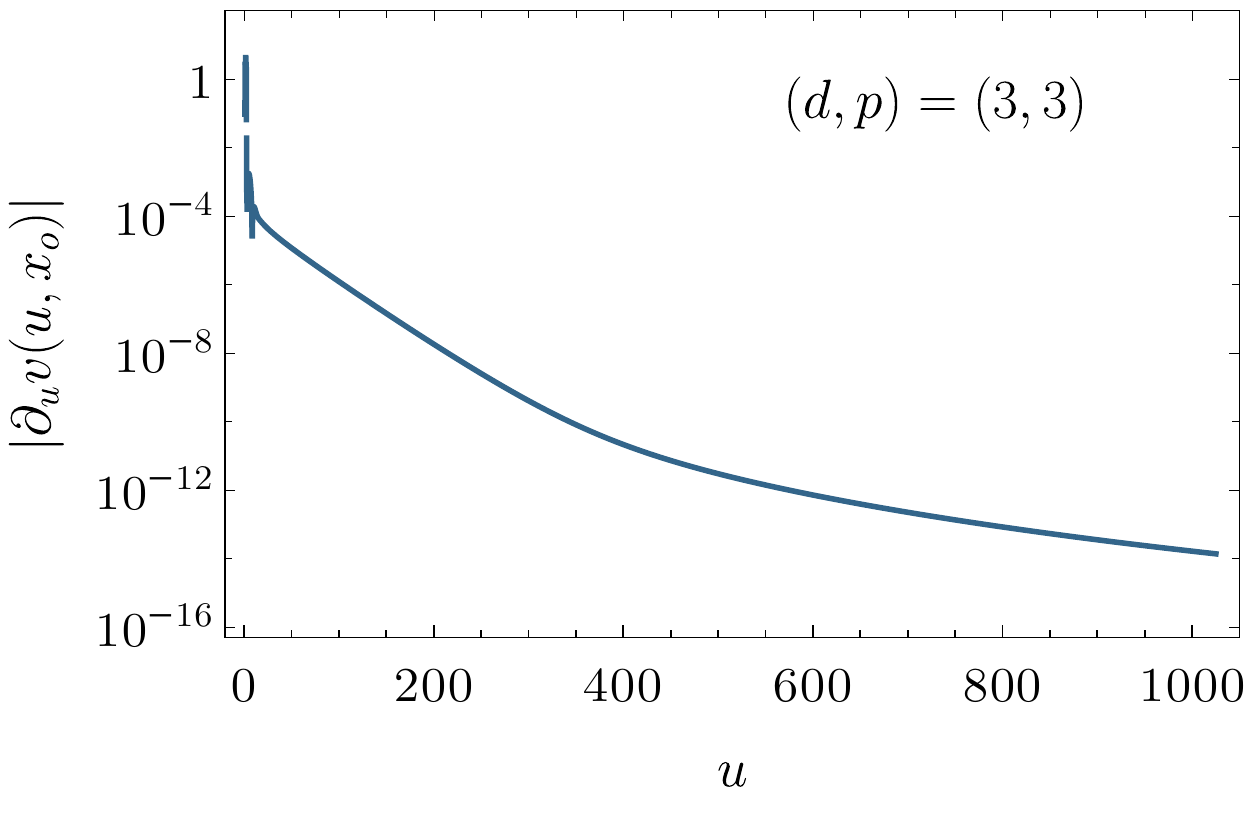}
  \hspace{1.5ex}
  \includegraphics[width=0.48\textwidth]{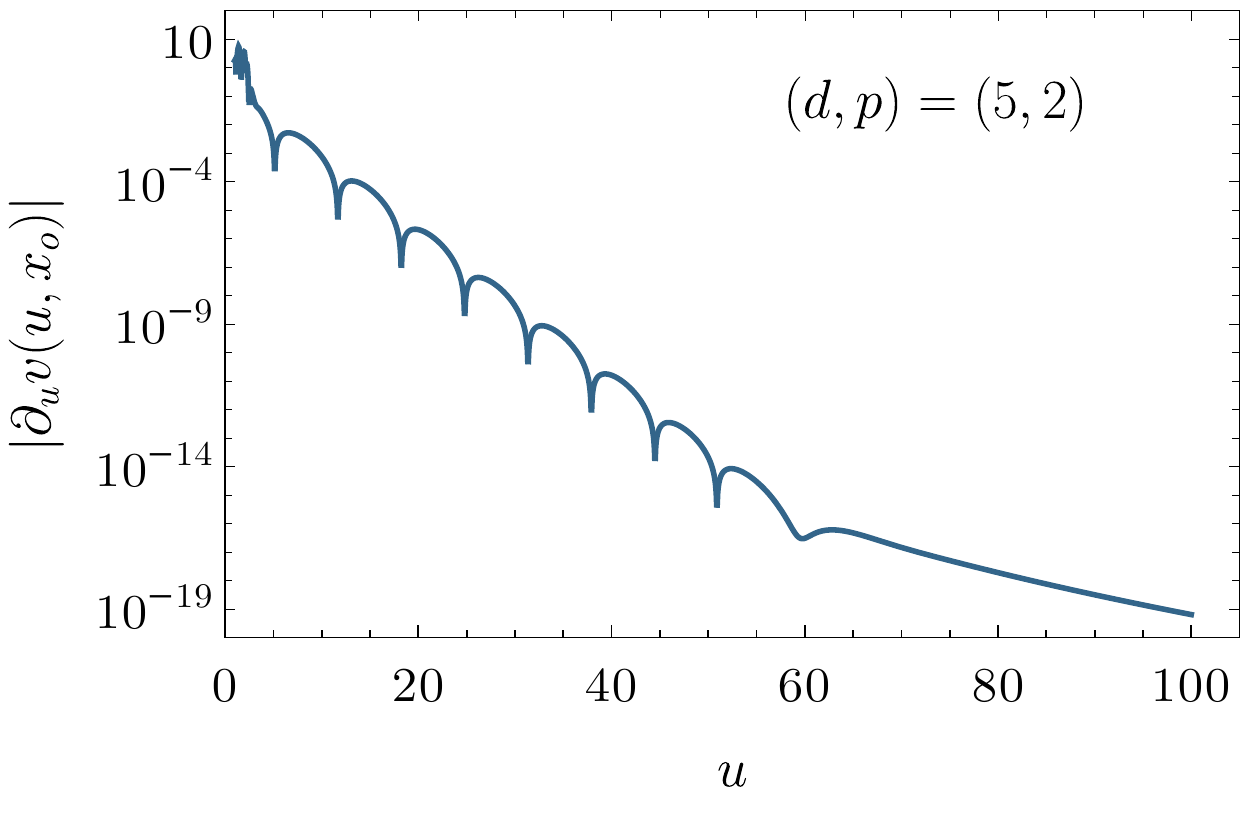}
  \scaption{Pointwise decay (at a sample interior point $x_{o}=5/9$)
    of the linear perturbation of $f_{0}$ as governed by
    Eq.~(\ref{eq:1}) for $(d,p)=(3,3)$ and $(d,p)=(5,2)$ for sample
    compactly supported initial data with the unstable mode
    removed. Numerical solution uses $N=1024$ Chebyshev polynomials,
    cf. Fig.~\ref{fig:2}.}
  \label{fig:a2}
\end{figure}

Interestingly enough, using a large enough number of polynomials we
were able (after removing the unstable mode from the initial data) to
see in the evolution the polynomial tail whose exponent is in
agreement with \cite{bcr}, cf.~Fig.~\ref{fig:a2}.



\end{document}